%% file: 0paper.tex
\documentclass[12pt]{article}
\usepackage{amsfonts}
\usepackage{amsmath}
\usepackage{epsfig}

\title{The Crisscross and the Cup:
  Two Short $3$-Twist Paper Moebius Bands}
\author{Brienne Elisabeth Brown and Richard Evan Schwartz \thanks{Supported by N.S.F. Grant DMS-2102802 and also a Mercator Fellowship.}}

\newtheorem{theorem}{Theorem}[section]

\newtheorem{conjecture}[theorem]{Conjecture}

\begin{document}
\maketitle

\begin{abstract}
  We introduce the crisscross and the cup, both
  of which are
  immersed $3$-twist
  polygonal paper Moebius
  band of aspect ratio $3$.  We explain
  why these two objects are limits of
  smooth embedded paper Moebius bands
  having knotted boundary.  We conjecture
  that any smooth embedded paper Moebius
  band with knotted boundary has
  aspect ratio greater than $3$.
  The crisscross is planar but the cup is not.
    \end{abstract}

\input{1intro}

\input{2criss}
\input{3cup}
\input{refs}

\end{document}

%% file: 1intro.tex
\section{Introduction}

Informally, a {\it paper Moebius band\/} is what
you get when you take a strip of paper, give it
an odd number of twists in space, then tape the ends together.
A formal definition is given e.g. in [{\bf S1\/}].
A related concept is that of a {\it folded ribbon knot\/}.
This is what you get when you take a paper strip,
fold it up so that the ends meet, and then press it into
the plane.
A formal definition of a folded ribbon knot is given e.g.
in [{\bf DL\/}].  The difference between
paper Moebius bands and folded ribbon knots is that
the former are smooth surfaces in space and the
latter are polygonal objects in the plane,
with some additional combinatorial data akin
to a knot crossing diagram.

We say that a paper Moebius band is
{\it multi-twisted\/} if the boundary loop
$\partial M_{\lambda}$ is a non-trivial knot.
One can make a similar definition for
folded ribbon knots.
The multi-tristed case corresponds to
giving the strip of paper at least $3$ twists.

This paper is a sequel to the papers
[{\bf S1\/}] and [{\bf S2\/}].
In [{\bf S1\/}], R.E.S. 
resolves the minimum aspect ratio
question for paper Moebius bands,
discussed in W. Wunderlich's 1962 paper
[{\bf W\/}] and then
formally conjectured by
B. Halpern and C. Weaver [{\bf HW\/}] in 1977.
(See [{\bf T\/}] for an English
translation of [{\bf W\/}].)
The so-called {\it triangular paper
  Moebius band\/}, whose aspect ratio is
$\lambda=\sqrt 3$, is the best one
can do.  This example has an
unknotted boundary, a perfect
equilateral triangle.  The triangular
Moebius band is
not quite a paper Moebius band because
it is neither smooth nor embedded,
but it has an interpretation as a
folded ribbon knot.

In [{\bf S2\/}] R.E.S.  establishes
a similar result for paper cylinders
having at least $2$ twists. In this case
the best one has aspect ratio $2$ and
folds $4$ times around a right-isosceles
triangle.  This result, and some soft
work involving smooth approximations,
resolves [{\bf DL\/}, Conjecture 39]
in the case $n=1$.

B.E.B.  subsequently
got interested in this
work and decided to find the minimum
aspect ratio for a multi-twisted paper Moebius band
by a physical
experiment: Make a loose $3$-twist
paper Moebius band and
then carefully pull it tight. She
found two surprising limits, which
she calls the {\it crosscross\/} and the
{\it cup\/}.  Both of these have
aspect ratio $3$, and they seem
to be optimal in the sense that
they minimize aspect ratio amonst
multi-twisted paper Moebius bands.

The crisscross has an interpretation
as a folded ribbon (un)knot.  The cup
is somehow closer to a basketweave.
In this note we will describe the
crosscross carefully and then sketch
how it may be approximated by smooth
paper Moebius bands having knotted
boundary.  We will describe the cup
in a bit less detail, and then
discuss heuristically some of the mathematical
implications of the cup.

One impressive thing about these
objects is that
their aspect ratio is much less than $3\sqrt 3$, which
is what one would get from a $3$-twist paper
Moebius band by wrapping $3$-times around
the triangular Moebius band in a fairly obvious
way.  Indeed, the famous hexaflexagon is just
such a $3$-fold wrapping.  For a while R.E.S.
thought that the hexaflexagon had minimum aspect ratio
amongst multi-twist paper Moebius bands.
This same speculation for folded ribbon knots
is the content of [{\bf DL\/}, Conjecture 22].
The crisscross and the cup
demolish the hexaflexagon
in terms of aspect ratio.  Since the crisscross is
also a folded ribbon knot, it gives a
counterexample to [{\bf DL\/}, Conjecture 22].

The existence of the crisscross (or the cup) and its
knotted approximations establishes the following result.
\begin{theorem}
  Let $\lambda^*$ denote the infimal
  aspect ratio of a multi-twisted
  paper Moebius band.  Then $\lambda^* \leq 3$.
\end{theorem}

Inspired by the conjecture made by Halpern and Weaver
in 1977 we make the following conjecture.  
\begin{conjecture}
  \label{min}
  A smooth embedded multi-twisted paper Moebius
  band has aspect ratio greater than $3$.
\end{conjecture}
The conjecture and the theorem would combine to show that
$\lambda^*=3$.  In the category of folded ribbon knots,
Conjecture \ref{min} is the same as
[{\bf DL\/}, Conjecture 22] except that
$3$ replaces $3 \sqrt 3$.

The existence of the cup and the crisscross together
rules out the possibility of a theorem like
[{\bf S1\/}, Triangular Limit Theorem].  Assuming
that Conjecture \ref{min} is true, a minimizing
sequence of examples does not have a unique limit.
It might converge to the crisscross and it might
converge to the cup.  How many possible limits
are there?  We think that probably there are just
these two limits, but we are prepared for a surprise.

We would like to thank Elizabeth Denne, Ben Halpern,
Curtis McMullen,
Charles Weaver, and Sergei Tabachnikov
for helpful discussions about topics related
to this paper.

\newpage

%% file: 2criss.tex
\section{The Crisscross}

\subsection{Basic Description}

The left side of Figure 1 shows the folding
pattern for the crisscross.  The middle picture
shows an intermediate stage of folding.
The right side shows the thing all folded up.
The strip of paper is aqua on one side and
magenta on the other.

\begin{center}
\resizebox{!}{2.8in}{\includegraphics{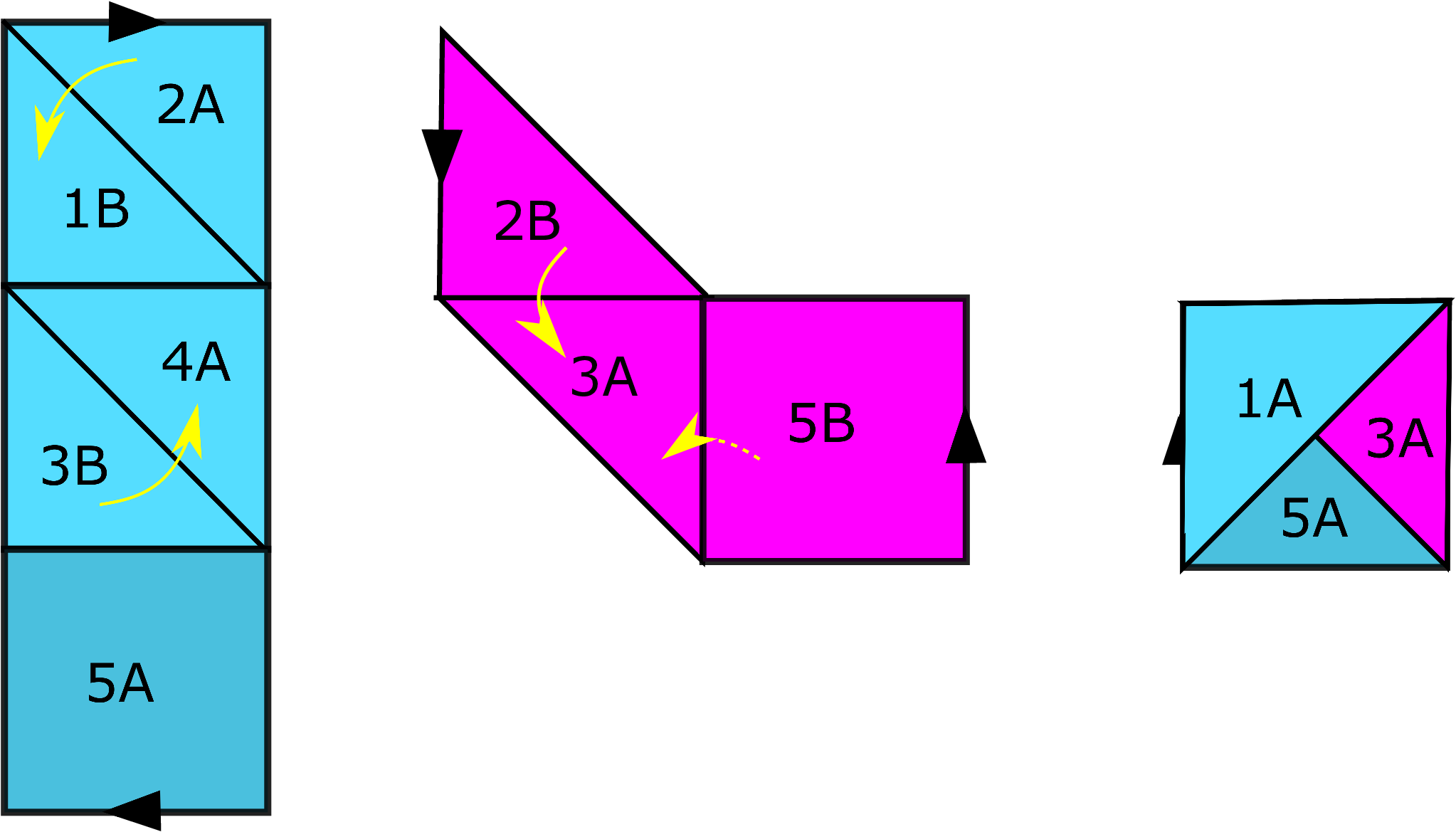}}
\newline
    {\bf Figure 1:\/} The crisscross folding pattern.
 \end{center}

The yellow arrows indicate the successive folds.
The solid yellow arrows indicate that the fold should
be made in a ``forward direction'' with the crease
receding away from the viewer. In the bottom
fold on the left, the bottom square goes along
for the ride. The one dotted
arrow in the middle indicates that the fold should
be made ``around the back''.  The lettering has
the following meaning:  Imagine that the crisscross
is siting on the table as on the right side of
Figure 1.  If you stick a pin though it, the pin
will encounter pieces $1,2,3,4,5$ in order.  The
$A$-faces are facing up and the $B$-faces are facing
down.  The crisscross is taped on the left hand
vertical side.

Figure 2 shows another view of the crosscross.
This time we are separating out the $5$ faces
and indicating how the edges are glued
together.  The thick red sides indicate
the boundary.  The oriented blue segments
piece together to make the midline of the
crisscross.

\begin{center}
\resizebox{!}{1in}{\includegraphics{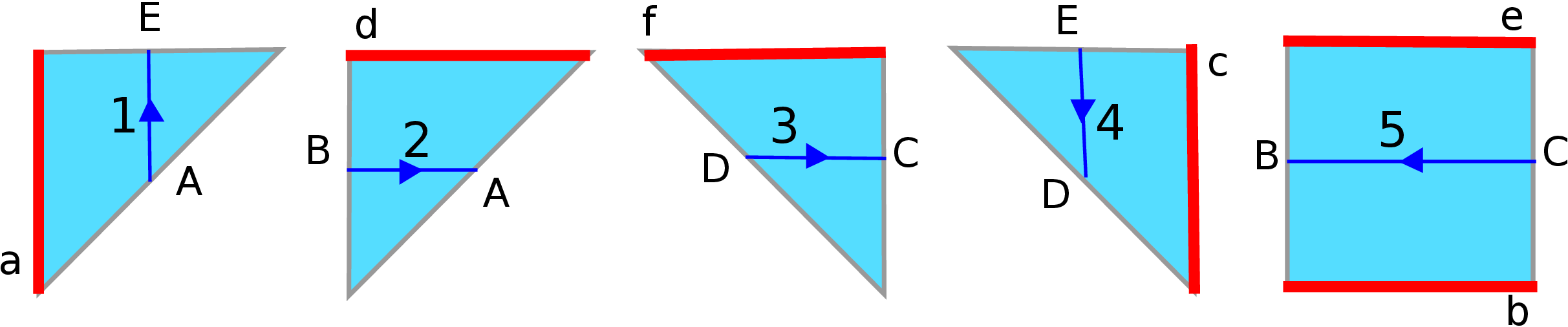}}
\newline
    {\bf Figure 2:\/} Crisscross gluing pattern
\end{center}

Figure 2 uses two systems of letters.  The letters
$A,B,C,D,E$ indicate the side pairings.  The
letters $a,b,c,d,e,f$ indicate the way the
boundary goes around the crisscross.  We
orient each of the sides labeled $a,b,c,d,e,f$
towards the letter.  Thus the tail of each edge is at the
unlettered vertex and the head is at the
lettered vertex.  The edge labeled
$a$ runs into the edge labeled $b$, with
the head of $a$ going to the tail of $b$.
And so on, all the way around.

\subsection{Smooth Approximation}

In this section we show that the crisscross
can be smoothly approximated as closely
as we like.  There are two ways to look
at this problem.  One way to look at it
is to observe that we made the crisscross
by making a $3$-twist Moebius band and pulling
it tight.  With a bit of faith in the
ability of mathematics to model physical
phenomena, one could imagine that the
intermediate stages are the
very approximations we seek.

The above approach is somewhat vague,
so here we describe a more precise method.
The method goes back at least
to the 1930 paper by M. Sadowski [{\bf Sa\/}].
(See [{\bf HF\/}] for an English translation.)
In Sadowski's paper the examples are made
by splicing together pieces of cylinders
with flat polygons, resulting in
$C^1$ examples with discontinuous
mean curvature.  A very similar
method, used by Halpern and Weaver [{\bf HW\/}],
uses more general smooth surfaces in place
of pieces of cylinders.
This argument is described informally
in [{\bf FT\/}].  We use the same argument,
more or less, in [{\bf S2\/}].
For convenience we repeat the argument
in [{\bf S2\/}] almost
{\it verbatim\/}.

Using
    smooth bump functions one can easily make a $U$-shaped
    curve.  This curve agrees with line parallel line
    segments at either end and then curves around to
    join these line segments.  Call this curve $U$.
    The product $U \times [a,b]$ is an isometrically
    embedded rectangle.  Next we take the polygons
    in Figure 2 and stack them on top of each other,
    separated by a very small distance.
    We now
    join the newly created edges of this stack
    of triangles by the $U$-shaped rectangles.
    
    If we make these $U$-shaped rectangles
    slightly thinner than the sides of the
    polygons, they will not overlap each other.
    The key observation here is that we never
    encounter $4$ edges laying directly above
    each other and having an interlaced gluing
    pattern.  Indeed, there are exactly $2$
    edges lying vertically above each other
    in each relevant plane.
    
    The new object will be an embedded smooth paper Moebius band
    whose boundary is not quite totally
    geodesic. To get a totally geodesic example, we just
    trim off the rough
    edges. This gives you a smooth embedded twisted
    cylinder with slightly larger aspect ratio.
    With this procedure you can make $\lambda$ as close
    to $3$ as you like.

    \subsection{Knottedness of the Boundary}

    We give $3$ ways to think about the knottedness of
    the boundary of the approximations to the crisscross.
    \newline
    \newline
    \noindent
    {\bf Method 1:\/}
    The informal method of approximating
    the crisscross by pulling tight a loose example, then
    the knotting is automatic.  The initial boundary is
    a trefoil knot and then when we pull tight we are
    moving the boundary by an isotopy and so it remains
    a trefoil knot.
    \newline
    \newline
    {\bf Method 2:\/}
    Another way to analyze the smooth approximations coming
    from the more formal construction is just to make one,
    then to tape some yarn along the edges, then
    detach the yarn, and then observe that the result
    is actually knotted.  We did this on the crisscross
    we built, which one can think of as a very close
    relative of the smooth approximations constructed above.

    Here we describe the boundary of an approximation to
    the crisscross based on the gluings in Figure 2.
    Figure 3 shows the path very nearly taken by the
    boundaries.  The numbers indicate the face containing
    the edges.  The little magenta segments indicate very
    nearly vertical segments which join a piece on one
    layer to a piece on another.  The crossings are
    dictated by the numbering.

\begin{center}
\resizebox{!}{2.7in}{\includegraphics{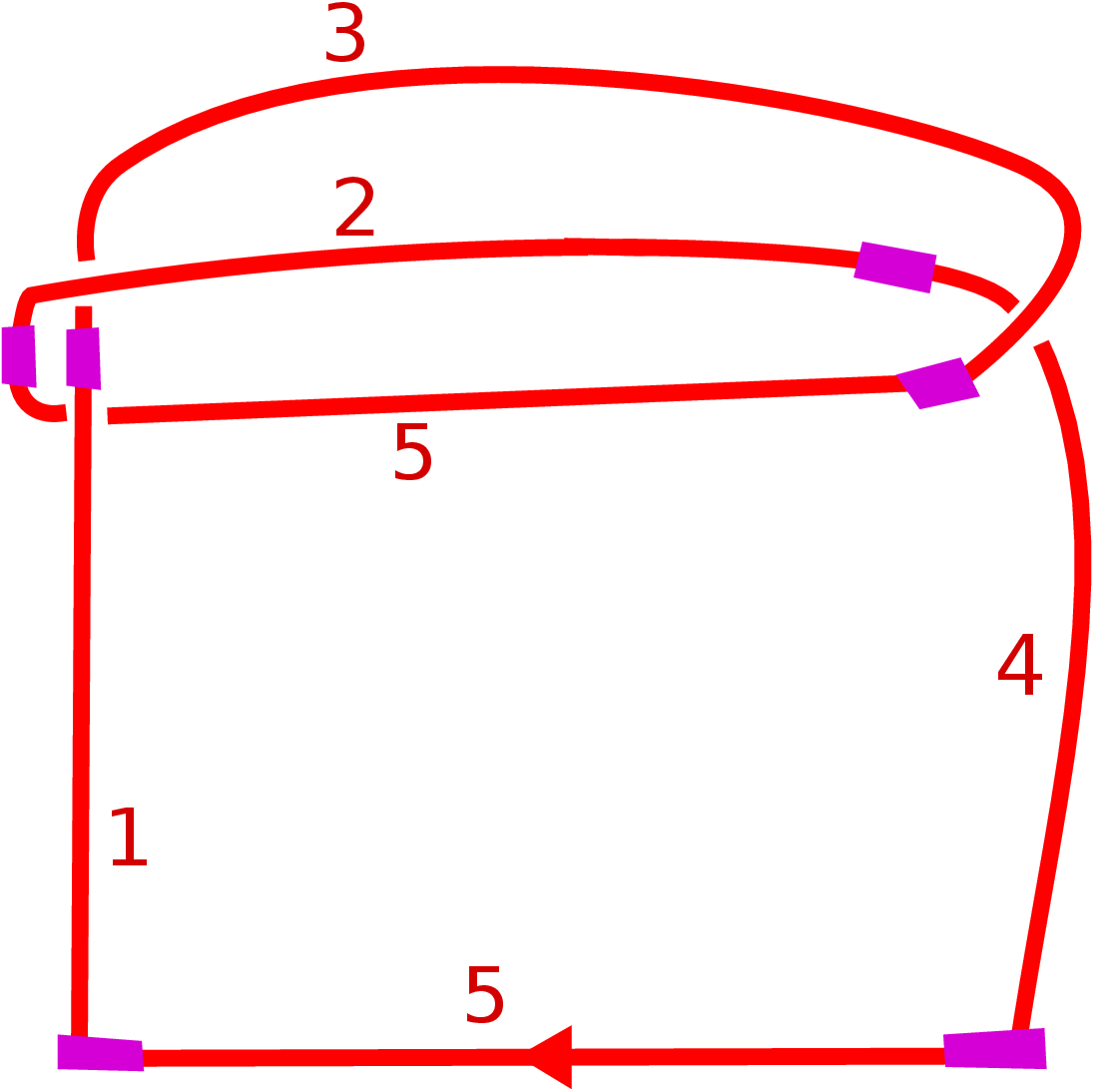}}
\newline
    {\bf Figure 3:\/} Knottedness of the boundary
\end{center}

In the actual crisscross the arcs which look horizontal
are horizontal and the arcs which look vertical are
vertical. The red  boundary makes a perfect square, but
with some backtracking.
\newline
\newline
    {\bf Method 3:\/}
    Our last method is completely algorithmic in
    the sense that it does not require any visualization
    or physical manipulation. It just requires an
    analysis of Figure 2.  In Figure 4 below we
    have copied down our knotted red loop, but
    for this method we do not need to know anything
    about this loop except the numbering of the
    strands.  The information
    inside the grey disks, which we got by
    some model-making in Method 2, is irrelevant
    for Method 3.
    
    The red loop in Figure 4 is as in Figure 3.  The
    blue loop is the midline of the crisscross.
    After orienting both the red and blue loops we
    can figure out which passes over which, at each
    crossing, using the numbering of the edges.
    Following this, we can assign a $(+)$ or a
    $(-)$ to each local crossing according to the
    rules given by the $16$ local models running
    around the outside of the picture.
    (One of the groups of $4$ comes from the Wikipedia
    page on linking number and the other groups are
    obtained from the first one by rotations.)
    Computing all the linking numbers and using the
    formula for linking number (total sign divided by
    $2$) we find that the linking number between the
    red/magenta and blue loops is $-3$.  This tells us
    that we have $3$-twist Moebius bands, and these
    have knotted boundaries.

\begin{center}
\resizebox{!}{5.5in}{\includegraphics{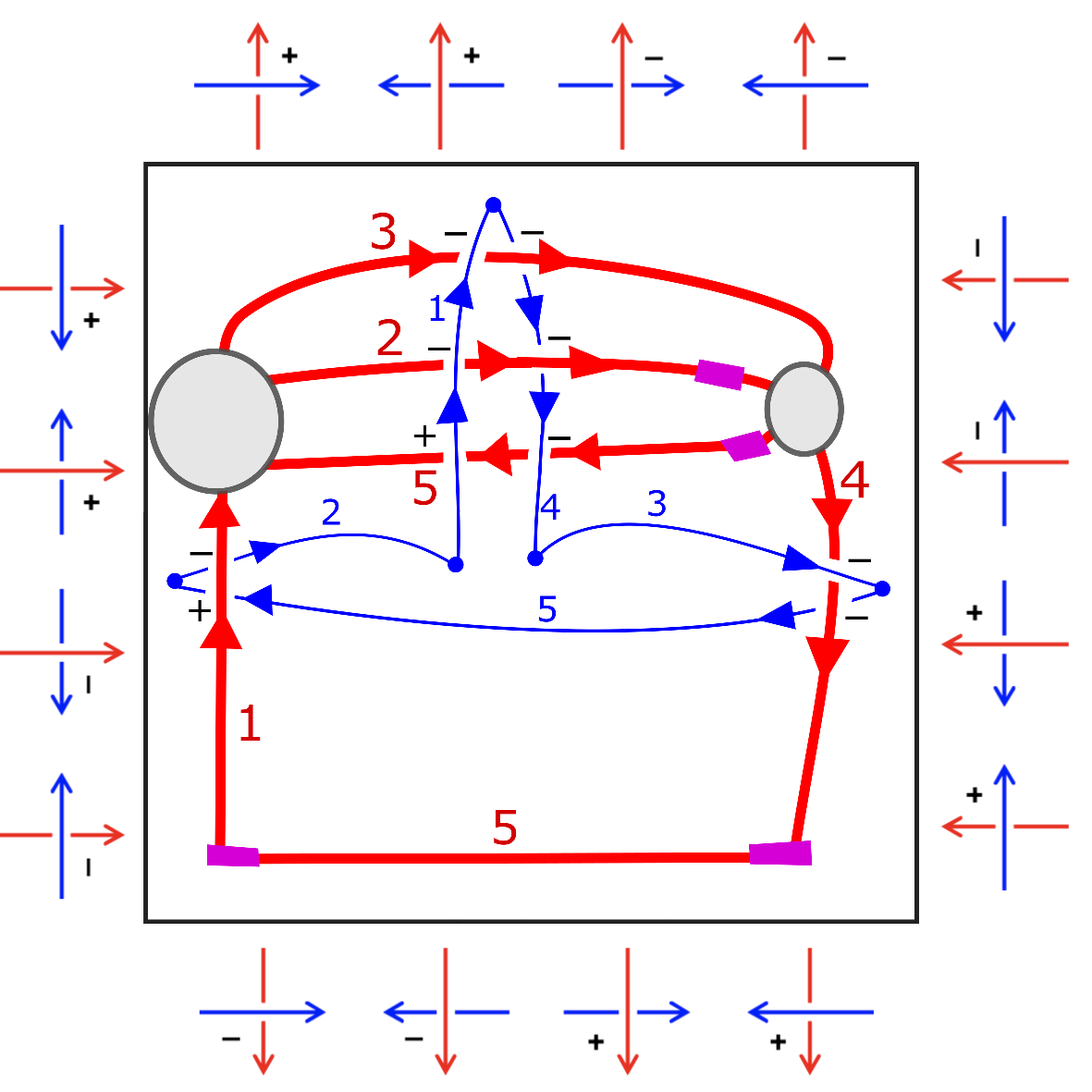}}
\newline
    {\bf Figure 4:\/} The linking of the boundary and the midline
\end{center}

Having more than one method for a relatively simple
example might seem like overkill, but it seems that
Method 3 is a nice way to algorithmically deduce
the twisting number of more general folded paper Moebius
bands based on gluing diagrams like Figure 2.
One should compare the method in
[{\bf DL\/}] for folded ribbon knots; it seems
very similar.

\newpage

%% file: 3cup.tex
\section{The Cup}

\subsection{Basic Description}

As we mentioned in the introduction,
B.E.B. has a much more
symmetric model for a $3$-twist
polygonal paper Moebius band of
aspect ratio $3$.  In this alternate
version, which we call
the {\it cup\/}, the image is a union of
$3$ right-isosceles triangles which make
$3$ faces of a tetrahedron whose fourth
face is an equilateral triangle.

\begin{center}
\resizebox{!}{3.2in}{\includegraphics{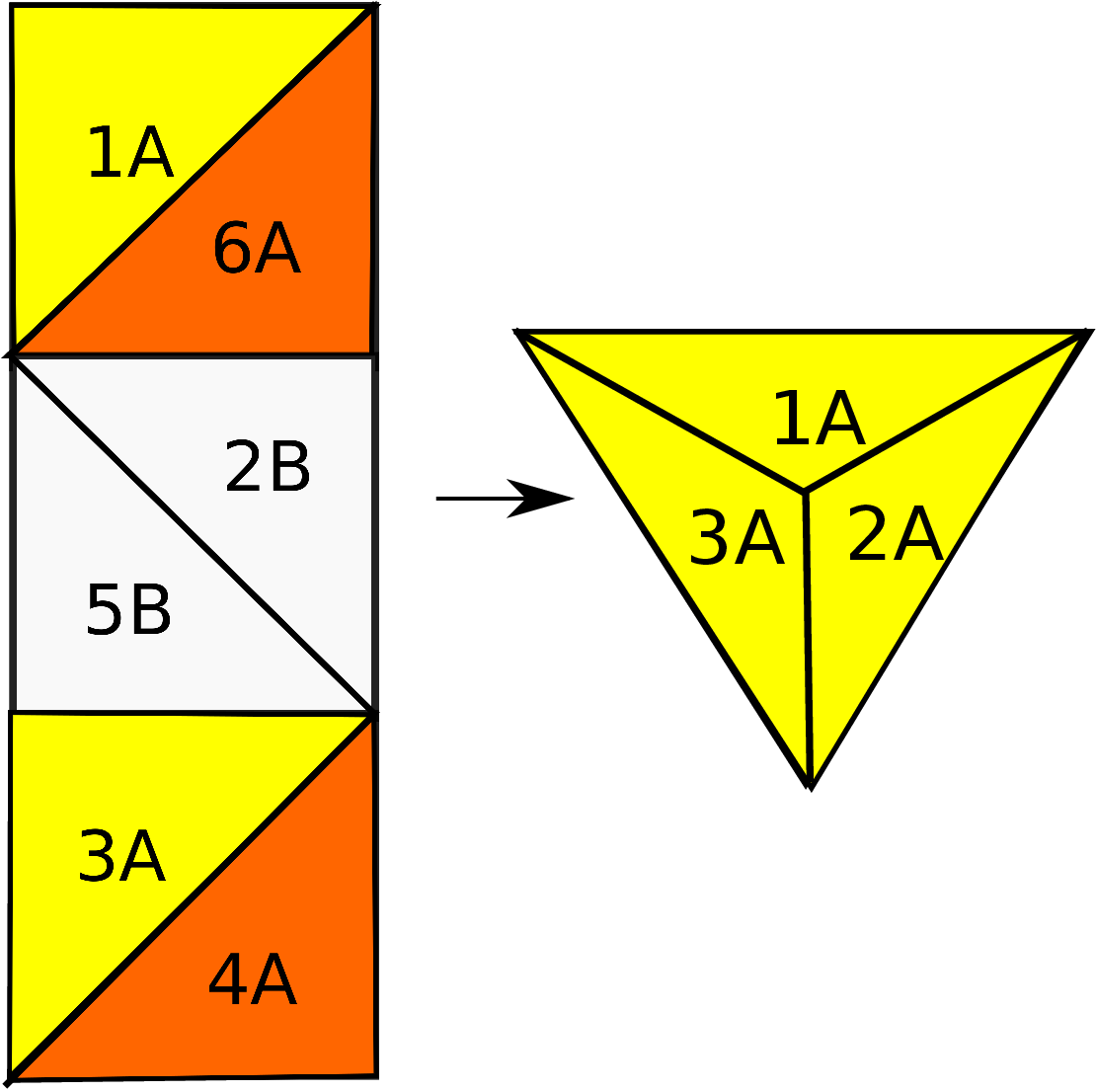}}
\newline
    {\bf Figure 5:\/} The folding pattern for the cup.
\end{center}

The left side of Figure 5 shows the
folding pattern for the cup. The right
side shows what you would see if you
were looking into the cup.  The lettering and
coloring is a bit different.  The
$A$-faces are colored orange and yellow and the
$B$-faces are colored white. When
this is folded up, the cup is
yellow on the ``inside'' and orange on
the ``outside''.  The white faces are
all pressed together and would not be
visible if you held the cup in your hands.
The faces $1A, 2A, 3A$ are on the ``inside'' of
the cup and the faces
$4A, 5A, 6A$ are on the ``outside''.
Again, the cup does not lie flat in
the plane and so the right side of
Figure 5 is a planar projection and
hence geometrically distorted.

We leave the details of the construction
to the reader who likes to cut out strips
of paper and fold them up.   Here is an
example we made from wax paper and
washi tape.

\begin{center}
\resizebox{!}{6.4in}{\includegraphics{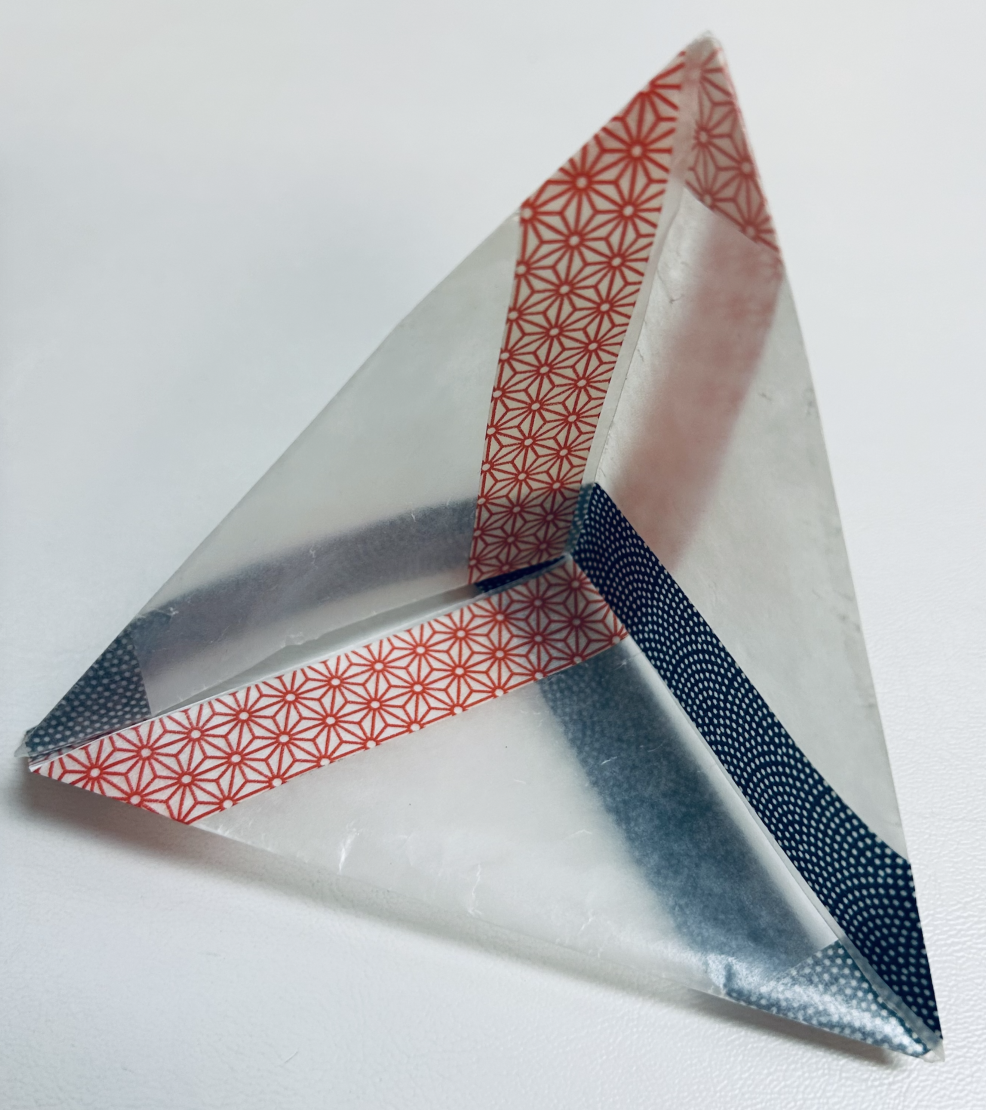}}
\newline
{\bf Figure 6:\/}  A cup made from wax paper and washi tape.
  \end{center}

\subsection{Smooth Approximations}

We approximate the cup by smooth paper
Moebius bands much in the same way as we
did for the crisscross.  We separate out
the $6$ triangular faces from each other
and then attach the appropriate
smooth folds.
The folding pattern is very interesting.  First of all,
there are $3$ folds which go around the rim of the cup.
Then there are $3$ more folds which make a $Y$-pattern
and meet at the central point of the cup.

It is interesting to observe how the boundary interacts with
the picture. The boundary consists of $6$ unit segments,
$3$ of which run along the inside of the cup in a $Y$-pattern
and $3$ of which run along the outside of the cup in a $Y$-pattern.
The whole picture has $3$-fold rotational symmetry, so if you
run your finger around either on the inside or the outside
of the cup you will encounter the boundary in a kind of
pinwheel fashion.

One can make the smooth approximations with
$3$-fold symmetry.  Figure 6 shows that the
boundary of such an approximation looks like.

\begin{center}
\resizebox{!}{4in}{\includegraphics{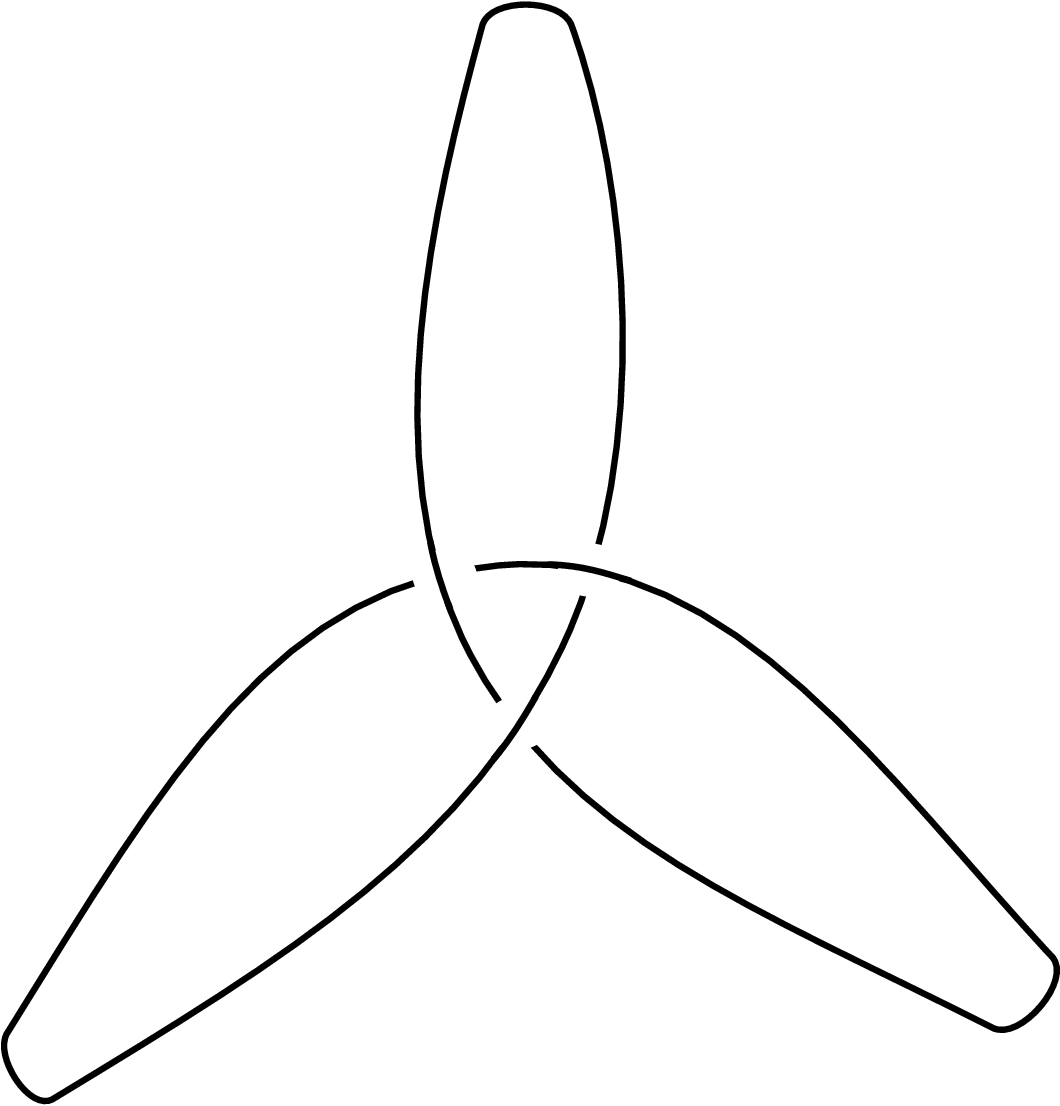}}
\newline
    {\bf Figure 7:\/} The boundary of a nearby smooth approximation.
  \end{center}

\subsection{Discussion}

The existence of the cup means that the kinds
of arguments made in [{\bf S1\/}] and [{\bf S2\/}],
by themselves, are unlikely to resolve Conjecture \ref{min}.
The idea in [{\bf S1\/}] is to first show that a
paper Moebius band has a pair of coplanar
and perpendicular bends, called a $T$-{\it pattern\/}.
After this, the idea is to consider the picture in
the plane of the $T$-pattern and establish an
optimization result.  The idea in [{\bf S2\/}]
is to find a nice planar projection of a twisted
paper cylinder and then do a similar kind of
optimization trick.  Both these approaches rely
on planar ideas.

It is worth pointing out that [{\bf S1\/}, Triangular Limit Theorem]
implies the existence of some $\eta>0$ such that any
knotted paper Moebius band has aspect ratio
at least $\sqrt 3 + \eta$.  The point is that any
paper Moebius band having aspect ratio very close
to $\sqrt 3$ is also very close to the triangular
Moebius band.  This would force the boundary of
the paper Moebius band to be unknotted, a
contradiction.   

More directly, a multi-twisted paper Moebius
band $\Omega$ has a $T$-pattern.  If the aspect ratio
of $\Omega$ is very near $\sqrt 3$ then the convex hull
$\nabla$ of this $T$-pattern is close to an
equilateral triangle of perimeter $2\sqrt 3$.
Since $\partial \Omega$  is knotted, the projection
of $\partial \Omega$ could not just follow along
$\partial \nabla$. 
This gives us the extra length and hence the better bound.

We did not try to find $\eta$, but
we can say that such a scheme will never prove
that $\sqrt 3 + \eta =3$. Here is the problem:  As a
byproduct of such an argument we would also prove that
any near minimizer for the aspect ratio is also nearly
planar.  But there are near minimizers which approximate
the cup, and these are far from planar.  This is a
contradiction.  

Another approach to Conjecture \ref{min}
would be to show that an arbitrary
paper Moebius band could be deformed, through
isometric embeddings,
into a folded ribbon knot.
Call such a paper Moebius band {\it flattenable\/}.
Conjecture
\ref{min} seems easier in the category of
folded ribbon knots, though probably still
quite hard.  If every paper Moebius band was
flattenable, the folded ribbon knot case
would imply the general case.

The cup is a rigid object.  Presumably, 
nearby smooth approximations are also not flattenable.
If they were, then we could probably extract a limit
and contradict the rigidity of the cup. (We have
not thought through the details of this.) In any
case, one can ask: When is
a paper Moebius band flattenable?  We don't know
any conditions which imply the answer one way
or the other.

\newpage

%% file: refs.tex
\section{References}

   [{\bf D\/}] E. Denne, {\it Folded Ribbon Knots in the Plane\/}, 
    The Encyclopedia of Knot Theory (ed. Colin Adams, Erica Flapan, Allison Henrich, Louis H. Kauffman, Lewis D. Ludwig, Sam Nelson)
    Chapter 88, CRC Press (2021)
    \vskip 8 pt
    \noindent
    [{\bf DL\/}] E. Denne, T. Larsen, {\it Linking number and folded ribbon unknots\/},
    Journal of Knot Theory and Its Ramifications, Vol. 32 No. 1 (2023)
    \vskip 8 pt
\noindent
[{\bf FT\/}], D. Fuchs, S. Tabachnikov, {\it Mathematical Omnibus: Thirty Lectures on Classic Mathematics\/}, AMS 2007
\vskip 7 pt
\noindent
[{\bf HF\/}], D.F. Hinz, E. Fried, {\it Translation of Michael Sadowsky’s paper 
‘An elementary proof for the existence of a developable MÖBIUS band and the attribution 
of the geometric problem to a variational problem’\/}. J. Elast. 119, 3–6 (2015)
\vskip 7 pt
\noindent
[{\bf HW\/}], B. Halpern and C. Weaver,
{\it Inverting a cylinder through isometric immersions and embeddings\/},
Trans. Am. Math. Soc {\bf 230\/}, pp 41--70 (1977)
\vskip 7 pt
\noindent
[{\bf Sa\/}], M. Sadowski, {\it Ein elementarer Beweis für die Existenz eines 
abwickelbaren MÖBIUSschen Bandes und die Zurückführung des geometrischen 
 Problems auf einVariationsproblem\/}. Sitzungsberichte der Preussischen Akad. der Wissenschaften, physikalisch-mathematische Klasse 22, 412–415.2 (1930)
\vskip 7 pt
\noindent
    [{\bf S1\/}] R. E. Schwartz, {\it The optimal paper Moebius Band\/},
    arXiv 2308.12641
\vskip 7 pt
\noindent
    [{\bf S2\/}] R. E. Schwartz, {\it The optimal twisted paper cylinder\/},
    arXiv 2309.14033
\vskip 7 pt
\noindent
[{\bf T\/}] Todres, {\it Translation of W. Wunderlich's On a Developable M\"obius band\/},
Journal of Elasticity {\bf 119\/} pp 23--34 (2015)
\vskip 7 pt
\noindent [{\bf W\/}] W. Wunderlich, {\it \"Uber ein abwickelbares M\"obiusband\/}, Monatshefte f\"ur Mathematik {\bf 66\/} pp 276--289 (1962)